\documentclass[a4paper,12pt]{article}

\usepackage{amsmath}
\usepackage{amsfonts}
\usepackage{amssymb}

\usepackage{hyperref}

\usepackage{color}

\usepackage[T1]{fontenc}
\usepackage{graphicx}

\usepackage{epsf}
\usepackage{psfrag}
\usepackage{subfigure}

\usepackage{fancyhdr}

\newtheorem{theorem}{Theorem}[section]
\newtheorem{lemma}[theorem]{Lemma}

\newtheorem{proposition}[theorem]{Proposition}
\newtheorem{remark}[theorem]{Remark}

\numberwithin{equation}{section}

\title{Invariant elliptic curves as attractors in the projective plane}
\author{Johan Taflin}

\begin{document}

\maketitle

\begin{abstract}
Let $f$ be a rational self-map of $\mathbb P^2$ which leaves invariant an elliptic curve $\mathcal C$ with strictly negative transverse Lyapunov exponent. We show that $\mathcal C$ is an attractor, i.e. it possesses a dense orbit and its basin has strictly positive measure.
\end{abstract}
\section{Introduction}

Let $f$ be a rational self-map of $\mathbb P^2$ of algebraic degree $d\geq2$ which leaves invariant an elliptic curve $\mathcal C$ (i.e. an algebraic curve of genus one). We assume that $\mathcal C$ does not contain indeterminacy points. In \cite{BDM}, Bonifant, Dabija and Milnor study such maps and give several examples. They associate to $f,$ a canonical ergodic measure $\mu_{\mathcal C}$, supported on $\mathcal C,$ which possesses a strictly positive Lyapunov exponent $\chi_1=(\log d)/2$ in the tangent direction of $\mathcal C.$ The transverse exponent corresponds to the second Lyapunov exponent $\chi_2$ of $\mu_\mathcal C,$ see section \ref{sec hyp} for the definition.

An invariant compact set $A=f(A)$ will be called \textit{an attractor} if $A$ possesses a dense orbit and if the basin of $A$ defined by
$$B(A)=\{x\in\mathbb P^2\,|\, d(f^n(x),A)\rightarrow 0\ \mbox{as}\ n\rightarrow\infty\},$$
has strictly positive Lebesgue measure. Here, $d(.,.)$ denotes the distance in $\mathbb P^2$ with respect to a fixed Riemanian metric. The purpose of this article is to establish the following theorem which was expected by Bonifant, Dabija and Milnor.
\begin{theorem}\label{th main}
Let $f$, $\mathcal C$ and $\mu_{\mathcal C}$ be as above. Assume that the transverse exponent $\chi_2$ of $\mu_{\mathcal C}$ is strictly negative. Then $\mathcal C$ is an attractor.
\end{theorem}
A sketch of proof is given in \cite{AKYY} using the absolute continuity of the stable foliation. Our strategy is to study the stable manifolds associated to $\mu_\mathcal C$, and then apply the following local result.
\begin{lemma}\label{le key lemma sur la positivite de l union des disque}
 Let $E$ be a subset of the unit disk $\Delta$ with strictly positive Lebesgue measure. Suppose $\{D_x\}_{x\in E}$ is a measurable family of disjoint holomorphic disks given by $\rho_x:\Delta\rightarrow\Delta^2$, transverse to $\{0\}\times\Delta$ and such that $\rho_x(0)=(0,x)$. Then the union $\cup_{x\in E}D_x$ has strictly positive Lebesgue measure in $\Delta^2$.
\end{lemma}
The proof of this lemma is based on holomorphic motions and quasi-conformal mappings.
Under the assumptions of Theorem \ref{th main}, $\mu_\mathcal C$ is a saddle measure, see \cite{deT}, \cite{Di} for the construction of such measures in a similar context and \cite{Si}, \cite{DS} for the basics on complex dynamics.

Recall that a rational self-map $f$ of $\mathbb P^2$ of algebraic degree $d$ is given in homogeneous coordinates $[z]=[z_0:z_1:z_2]$, by $f[z]=[F_0(z):F_1(z):F_2(z)]$ where $F_0,F_1,F_2$ are three homogeneous polynomials in $z$ of degree $d$ with no common factor. In the sequel, we always assume that $d\geq 2$. The common zeros in $\mathbb P^2$ of $F_0,$ $F_1,$ and $F_2$ form the indeterminacy set $I(f)$ which is finite.
Let $\mathcal C\subset\mathbb P^2$ be an elliptic curve. Then, there exists a lattice $\Gamma$ of $\mathbb C$ and a desingularization
$$\Psi: \mathbb C/\Gamma\rightarrow\mathbb P^2,$$
with $\Psi(\mathbb C/\Gamma)=\mathcal C.$ Moreover, if $S$ denotes the singular locus of $\mathcal C,$ the map
$$\Psi: (\mathbb C/\Gamma)\setminus \Psi^{-1}(S)\rightarrow \mathcal C\setminus S$$
is a biholomorphism.

We say that $\mathcal C$ is $f$-\textit{invariant} if $\mathcal C\cap I(f)=\varnothing$ and $f(\mathcal C)=\mathcal C.$ In this case, the restriction $f_{|\mathcal C}$ lifts to a holomorphic self-map $\widetilde f$ of $\mathbb C/\Gamma.$
Even if $\mathcal C$ is singular, $f$ inherits several properties of $\widetilde f.$ Like all holomorphic self-maps of $\mathbb C/\Gamma$, $\widetilde f$ is necessarily of the form $t\mapsto at+b$ and leaves invariant the normalized Lebesgue measure $\widetilde\mu_{\mathcal C}$ on $\mathbb C/\Gamma.$ So, the topological degree of $\widetilde f$, i.e. the number of points in a fiber, is equal to $|a|^2$. It is not difficult to check that this degree is equal to $d$, see \cite{BD}. Therefore, $|a|^2=d\geq2$.
Then, by a classical theorem on ergodicity on compact abelian groups, $\widetilde\mu_{\mathcal C}$ is $\widetilde f$-ergodic, i.e. is extremal in the cone of invariant positive measures. Its push-forward $\mu_{\mathcal C}$ is an $f$-ergodic measure supported on $\mathcal C.$ Moreover, generic orbits of $f_{|\mathcal C}$ are dense in $\mathcal C.$ On the other hand, $f_{|\mathcal C}$ inherits the repulsive behavior of $\widetilde f$ and $\mu_\mathcal C$ possesses a strictly positive Lyapunov exponent equal to $\chi_1=\log|a|=(\log d)/2$ in the tangent direction of $\mathcal C.$ 
By Oseledec's theorem, see Section \ref{sec hyp} below, we have 
$$\chi_1+\chi_2=\dfrac{1}{2}\langle\mu_\mathcal C, \log(\mbox{Jac}(f))\rangle$$
where $\mbox{Jac}(f)$ denotes the Jacobian of $f$ with respect of the Lebesgue measure of $\mathbb P^2$. So, the hypothesis in Theorem \ref{th main} is equivalent to
$$\langle \mu_\mathcal C,\log(\mbox{Jac}(f))\rangle <\log d.$$
Some examples in \cite{BDM} satisfy this condition and give the first attractors in $\mathbb P^2$ with non-open basins.
\vspace{1ex}
\\
\textbf{Acknowledgements.} I would like to thank Tien-Cuong Dinh for drawing my attention to the subjet and for his invaluable help.

\section{Holomorphic motion}
We briefly introduce the notion of holomorphic motion. For a more complete account cf. \cite{GJW}. For $r>0,$ we denote by $\Delta_r$ the disk centered at the origin in $\mathbb C$ with radius $r$. If $E$ is a subset of $\mathbb P^1$, a holomorphic motion of $E$ parametrized by $\Delta$ is a map
$$h:\Delta\times E\rightarrow \mathbb P^1$$
such that:
\begin{enumerate}
 \item[i)] $h(0,z)=z$ for all $z\in E$,
 \item[ii)] $\forall \, c\in\Delta$, $z\mapsto h(c,z)$ is injective,
 \item[iii)] $\forall \,z\in E$, $c\mapsto h(c,z)$ is holomorphic on $\Delta$.
\end{enumerate}
By the works of Ma\~né, Sad, Sullivan, Thurston and Slodkowski (see \cite{MSS}, \cite{ST} and \cite{Slo}), any holomorphic motion $h$ of $E$ can be extended to a holomorphic motion $\widetilde h$ of $\mathbb P^1.$ Furthermore, even if no continuity in $z$ is assumed, $\widetilde h$ is continuous on $\Delta\times\mathbb P^1.$ More precisely, we have the following result.
\begin{theorem}\label{th holo mot}
Let $h$ be a holomorphic motion of a set $E\subset\mathbb P^1$ parametrized by $\Delta$. Then there is a continuous holomorphic motion $\widetilde h:\Delta\times\mathbb P^1\rightarrow\mathbb P^1$ which extends $h$. Moreover, for any fixed $c\in\Delta$, $\widetilde h(c,.):\mathbb P^1\rightarrow\mathbb P^1$ is a quasi-conformal homeomorphism.
\end{theorem}
We refer to \cite{Ahl} for quasi-conformal mappings. The following property is crucial in our proof.
\begin{proposition}\label{prop qc mapping et absolute continuity}
A quasi-conformal mapping sends sets of Lebesgue measure $0$ to sets of Lebesgue measure $0$.
\end{proposition}
We shall need the following Lemma in the proof of Lemma \ref{le key lemma sur la positivite de l union des disque}.
\begin{lemma}\label{le mesure de l'union}
 Let $h$ be a holomorphic motion of a Borel set $E\subset\mathbb P^1$ of strictly positive measure. Then $\cup_{c\in\Delta}\{c\}\times h(c,E)$ has strictly positive measure in $\Delta\times\mathbb P^1$.
\end{lemma}
\textit{Proof}. By Theorem \ref{th holo mot}, $h$ can be extended to a holomorphic motion $\widetilde h: \Delta\times \mathbb P^1\rightarrow\mathbb P^1$ such that, for any fixed $c\in\Delta$, $\widetilde h(c,.):\mathbb P^1\rightarrow\mathbb P^1$ is a quasi-conformal homeomorphism. So, by Proposition \ref{prop qc mapping et absolute continuity} $\mbox{Leb}(h(c,E))>0$ and Fubini's theorem implies that $\cup_{c\in\Delta}\{c\}\times h(c,E)$ has strictly positive measure.
\hfill $\Box$

\vspace{0.5cm}
\noindent\textit{Proof of Lemma \ref{le key lemma sur la positivite de l union des disque}}. From the family of disks, we will construct a holomorphic motion. Denote by $\pi_1$ and $\pi_2$ the canonical projections of $\Delta^2$. Let $x\in E$. Since $D_x$ is transverse to $\{0\}\times\Delta$, there exists $r(x)>0$ such that $\rho^1_x=\pi_1\circ\rho_x$ is a biholomorphism between a neighbourhood of $0$ and $\Delta_{r(x)}$. The measurability of $\{D_x\}_{x\in E}$ implies that $x\mapsto r(x)$ is also measurable. As $\mbox{Leb}(E)>0$, there exists $a>0$ and a subset $E_a$ of $E$ such that $\mbox{Leb}(E_a)>0$ and $r(x)>a$ for each point $x\in E_a.$ Define
\begin{align*}
 h:\Delta_a&\times E_a\rightarrow \Delta\\
(c&,z)\mapsto \rho_{z}^2\circ(\rho_{z}^1)^{-1}(c),
\end{align*}
where $\rho^2_x=\pi_2\circ\rho_x$. By construction, $h$ is well defined and $c\mapsto h(c,z)$ is holomorphic on $\Delta_a$. Since the disks are pair-wise disjoint, the map $z\mapsto h(c,z)$ is injective for each $c\in\Delta_a$. Therefore, $h$ is a holomorphic motion of $E_a$ parametrized by $\Delta_a$ and by Lemma \ref{le mesure de l'union} $\cup_{c\in\Delta_a}\{c\}\times h(c,E_a)\subset \cup_{x\in E}D_x$ has strictly positive Lebesgue measure in $\Delta^2$.
\hfill $\Box$
\section{Hyperbolic dynamics}\label{sec hyp}
Suppose that $g$ is a holomorphic self-map of a complex manifold $M$ of dimension $m$.
The following Oseledec's multiplicative ergodic theorem (cf. \cite{KH} and \cite{Wa}) gives information on the growth rate of $\|D_xg^n(v)\|,$ $v\in T_xM$ as $n\rightarrow +\infty.$ Here, $D_xg^n$ denotes the differential of $g^n$ at $x.$ Oseledec's theorem holds also when $g$ is only defined in a neighbourhood of $\mbox{supp}(\nu).$
\begin{theorem}\label{th ose}
 Let $g$ be as above and let $\nu$ be an ergodic probability with compact support in $M$. Assume that $\log^+ \mbox{Jac}(g)$ is in $L^1(\nu)$. Then there exist integers $k,$  $m_1,\ldots,m_k,$ real numbers $\lambda_1>\cdots>\lambda_k$ ($\lambda_k$ may be $-\infty$) and a subset $\Lambda\subset M$ such that $g(\Lambda)=\Lambda$, $\nu(\Lambda)=1$ and for each $x\in\Lambda,$ $T_xM$ admits a measurable splitting
$$T_xM=\bigoplus_{i=1}^kE^i_x$$
such that $\dim_{\mathbb C}(E_x^i)=m_i,$ $D_xg(E^i_x)\subset E_{g(x)}^i$ and
$$\lim_{n\rightarrow +\infty}\dfrac{1}{n}\log\|D_xg^n(v)\|=\lambda_i$$
locally uniformly on $v\in E_x^i\setminus\{0\}.$ Moreover, for $S\subset N:=\{1,...,k\}$ and $E_x^S=\oplus_{i\in S}E_x^i,$ the angle between $E^S_{g^n(x)}$ and $E^{N\setminus S}_{g^n(x)}$ satisfies
$$\lim_{n\rightarrow +\infty}\dfrac{1}{n}\log\sin|\angle(E^S_{g^n(x)},E^{N\setminus S}_{g^n(x)})|=0.$$
\end{theorem}
\vspace{2ex}
The constants $\lambda_i$ are the \textit{Lyapunov exponents} of $g$ with respect to $\nu.$ It is not difficult to deduce that
$$2\sum_{i=1}^km_i\lambda_i=\int \log \mbox{Jac}(g)d\nu.$$
If all Lyapunov exponents are non-zero, we say that $\nu$ is \textit{hyperbolic}. In this case, let $\lambda>0$ such that $\lambda<|\lambda_i|$ for all $1\leq i\leq k$ and let
$$E_x^s=\bigoplus_{\lambda_i<0}E^i_x,\ E_x^u=\bigoplus_{\lambda_i>0}E^i_x.$$
Then, for each point $x$ in $\Lambda$ and $\delta>0$ we define the \textit{stable manifolds} at $x$ by
$$W_\delta^s(x)=\{y\in M\,|\, d(g^n(x),g^n(y))<\delta e^{-\lambda n} \ \ \forall n\geq 0\}.$$
From Pesin's theory, we have the following fundamental result, see \cite{BP}, \cite{PS}, \cite{RS} and \cite{Pol} for more details.
\begin{theorem}\label{th sta mani}
There exists a strictly positive measurable function $\delta$ on $\Lambda$ such that if $x\in\Lambda$ then
\begin{itemize}
 \item[i)] $W_{\delta(x)}^s(x)$ is an immersed manifold in $M$,
 \item[ii)] $T_xW_{\delta(x)}^s(x)=E^s_x$,
 \item[iii)] $W_{\delta(x)}^s(x)$ depends measurably of $x$.
\end{itemize}
\end{theorem}

\section{Basin of an attracting curve}

Since the support of $\mu_{\mathcal C}$ does not contain indeterminacy points, we have $\log^+ \mbox{Jac}(f)\in L^1(\mu_{\mathcal C}).$ We assume that $\mu_{\mathcal C}$ has a strictly negative transverse exponent $\chi_2.$ Then, there exists a hyperbolic set $\Lambda\subset\mathcal C$ such that $\mu_{\mathcal C}(\Lambda)=1$ and $E_x^u=T_x\mathcal C$ for all $x\in\Lambda.$

The first step to apply Lemma \ref{le key lemma sur la positivite de l union des disque} is to find, for some $p\in\Lambda,$ an open neighbourhood where the stable manifolds are pair-wise disjoint. To this end, we prove that the restriction $f_{|\mathcal C}$ inherits the repulsive behavior of $\widetilde f.$ Recall that $d(.\,,.)$ is the distance on $\mathbb P^2$ and denote by $\widetilde d(.\,,.)$ the standard distance on $\mathbb C/\Gamma.$
\begin{lemma}\label{le repulsif}
There is a constant $\beta>0$ such that for each $p\in\mathcal C\setminus S$ we can find $\alpha>0$ with the property that, if $x,y\in \mathcal C$ are distinct points in the ball $B(p,\alpha)$ of radius $\alpha$ centered at $p$, then $d(f^n(x),f^n(y))>\beta$ for some $n\geq 0.$
\end{lemma}
\textit{Proof}. As $\Psi$ is one-to-one except on finitely many points, we can find a finite open covering $\{U_j\}_{j\in J}$ of $\mathbb C/\Gamma$ such that $\Psi$ is injective on each $\overline{U_j}.$ Let $z_1,z_2\in\mathbb C/\Gamma.$ We denote by $\epsilon>0$ a Lebesgue number of this covering, i.e. if $\widetilde d(z_1,z_2)<\epsilon$ then, there exists $j\in J$ such that $z_1$ and $z_2$ are in $U_j.$ Recall that one can choose $r>0$ such that if $\widetilde d(z_1,z_2)<r$ then $\widetilde d(\widetilde f(z_1),\widetilde f(z_2))=|a|\widetilde d(z_1,z_2).$ We can assume that $\epsilon<r.$ Let $\widetilde\alpha>0$ such that $2|a|\widetilde\alpha\leq\epsilon.$ If $0<\widetilde d(z_1,z_2)<\widetilde\alpha$ then there exists $n\geq 0$ such that
$$\widetilde\alpha<\widetilde d(\widetilde f^n(z_1),\widetilde f^n(z_2))\leq |a|\widetilde\alpha.$$
Therefore, we can find $j\in J$ such that $\widetilde f^n(z_1)$ and $\widetilde f^n(z_2)$ are in $U_j.$ So $$d(f^n(\Psi(z_1)),f^n(\Psi(z_2)))>\beta,$$ where
$$\beta=\min_{j\in J}\inf_{\substack{x_1,x_2\in U_j\\\widetilde d(x_1,x_2)>\widetilde\alpha}}d(\Psi(x_1),\Psi(x_2))>0.$$
Finally, for each $p\in\mathcal C\setminus S$ we can choose $\alpha>0$ such that if $x,y\in B(p,\alpha)\cap\mathcal C,$ there are preimages $\widetilde x$ and $\widetilde y$ of $x$ and $y$ by $\Psi$ which satisfy $\widetilde d(\widetilde x,\widetilde y)<\widetilde\alpha.$ Then, $d(f^n(x),f^n(y))>\beta$ for some $n\geq 0.$
\hfill $\Box$
\begin{lemma}\label{le var disjointes}
Let $p\in\Lambda.$ There exist $\delta_0>0$ and an open neighbourhood $U$ of $p$ such that if $\delta<\delta_0,$ $x,y\in U\cap\Lambda,$ $x\neq y$ then $W_{\delta}^s(x)\cap W_{\delta}^s(y)=\varnothing.$
\end{lemma}
\textit{Proof}. By Lemma \ref{le repulsif}, we can choose for $U$ the ball of radius $\alpha$ centered at $p$ and $\delta_0\leq\beta/2.$ If there exist $x,y\in U\cap\Lambda,$ $x\neq y,$ with $W_{\delta}^s(x)\cap W_{\delta}^s(y)\neq\varnothing,$ then for $z\in W_{\delta}^s(x)\cap W_{\delta}^s(y),$
$$d(f^n(x),f^n(y))\leq d(f^n(x),f^n(z))+d(f^n(z),f^n(y))\leq 2\delta e^{-\lambda n}\leq 2\delta,$$
for every $n,$ which contradicts Lemma \ref{le repulsif}.
\hfill $\Box$
\vspace{2ex}\\
\textit{Proof of Theorem \ref{th main}}. Let $p\in \Lambda$ be a regular point of $\mathcal C$. Choosing suitable local coordinates at $p$, we can assume that $p\in \Delta^2$ and $\mathcal C\cap\Delta^2=\{0\}\times\Delta$. Let $x\in \Lambda\cap \Delta^2$. By the stable manifold theorem, there exists $\delta(x)>0$ such that $W_{\delta(x)}^s(x)\cap \Delta^2$ is an immersed manifold in $\Delta^2$. So, there exists a measurable family of embedded holomorphic disks $\rho_x:\Delta\rightarrow U$ with $\rho_x(0)=x$ and $\rho_x(\Delta)\subset W_{\delta(x)}^s(x)\cap U.$

First, by Lemma \ref{le var disjointes}, possibly after replacing $\Delta^2$ by a smaller polydisk, we can choose $\delta(x)<\delta_0$ for all $x\in\Lambda\cap \Delta^2$. The stable manifolds $W_{\delta(x)}^s(x)$ are then pair-wise disjoint.

Since $W_{\delta(x)}^s(x)$ is tangent to $E_x^s$ in $x$, the family of disks is transverse to $\{0\}\times\Delta$. Then, by Lemma \ref{le key lemma sur la positivite de l union des disque} the union of stable manifolds, which is included in the basin of $\mathcal C$, has strictly positive measure.
\hfill $\Box$

\begin{remark}
 By Hurwitz's formula, if $\mathcal C$ is an invariant curve then $\mathcal C$ is rational or elliptic. If $\mathcal C$ is rational and $f_{|\mathcal C}$ is a Lattès map, i.e a map which is semi-conjugated to an endomorphism of a torus, then its equilibrium measure is absolutely continuous with the respect to the Lebesgue measure. We obtain in the same way that its basin has strictly positive measure.
\end{remark}

\noindent
J. Taflin, UPMC Univ Paris 06, UMR 7586, Institut de
Math{\'e}matiques de Jussieu, F-75005 Paris, France. {\tt  taflin@math.jussieu.fr} 

\end{document}